\documentclass[3p,times]{elsarticle}

\usepackage{amsmath}

\usepackage{amssymb}

\usepackage[english]{babel}
\begin{document}
\large
\begin{frontmatter}
\large
		
\title{Theory of Generalized Trigonometric Functions: from Laguerre to Airy Forms}

\author[Enea]{G. Dattoli \corref{cor}}
\ead{giuseppe.dattoli@enea.it}
		
\author[Enea,Unict]{S. Licciardi }
\ead{silvia.licciardi@dmi.unict.it}
		
\author[Unict]{R.M. Pidatella}
\ead{rosa@dmi.unict.it}
		
\address[Enea]{ENEA - Frascati Research Center, Via Enrico Fermi 45, 00044, Frascati, Rome, Italy}
\cortext[cor]{Corresponding author}
\address[Unict]{Dep. of Mathematics and Computer Science, University of Catania, Viale A. Doria 6, 95125, Catania, Italy}

\begin{abstract}
We develop a new point of view to introduce families of functions, which can be identified as generalization of the ordinary trigonometric or hyperbolic functions. They are defined using a procedure based on umbral methods, inspired to the Bessel Calculus of Bochner, Cholewinsky and Haimo. We propose further extensions of the method and of the relevant concepts as well and obtain new families of integral transform allowing the framing of the previous concepts within the context of generalized Borel tranforms.
\end{abstract}
		
		
\end{frontmatter}
	
\section{Introduction}

\large

\noindent The cylindrical Bessel can be considered generalizations of the trigonometric functions, while the associated modified forms are an extension of the relevant hyperbolic counterparts \cite{L.C.Andrews}.\\
 Such an academic identification is non-particularly deep and might be useful for pedagogical reasons or as a guiding element to study their properties as e.g. those relevant to the asymptotic forms. We must however underline that Bessel and trigonometric/hyperbolic functions share some resemblances only, but they do not display any  full correspondence.\\
 The search for functions which are "true" generalizations of the trigonometric (t-) or hyperbolic (h-) forms is however recurrent in the mathematical literature.  The attempts in this direction can be ascribed to different strategies, roughly speaking the geometrical \cite{E.Ferrari} and the analytical \cite{D.E.Edmunds} point of views.\\
 The first is based on definitions extending to higher powers the Pythagorean identity of ordinary trigonometric functions, such a program identifies new trigonometries, with their own geometrical interpretation on elliptic curves and with different numbers playing the role of $\pi$ \cite{E.Ferrari}.\\
 The second invokes the analogy with series expansions, differential equations and the theory of special functions. \\
 
 The generalized $t\!-\!h$ functions, defined within these two contexts, are different; in particular those belonging to the geometric strategy can be recognized as Elliptic functions, including Jacobi and Weierstrass forms.  In this article we develop a systematic procedure within the framework of the analytical point of view.\\
 
 We look for "true" generalizations, in the sense that the functions we define allow a one to one mapping onto the properties of the elementary $t\!-\!h$ functions, like addition or duplication theorems. To this aim we exploit the insightful point of view offered by the recent understanding of Bessel functions as umbral manifestation of Gauss or of exponential functions \cite{Gorska}. These conceptual tools, as well as the ideas developed by Cholewinsky and Reneke in ref. \cite{Cholewinski}, provide the elements underlying the formalism of this paper, aimed at exploring in depth the identification of trigonometric functions associated with  Bessel functions, by getting the proper algebraic environment to establish the relevant properties.\\

 Our starting point is the equation
 
\begin{equation}\begin{split}\label{GrindEQ__1_} 
& {}_{l} \partial _{x} F(x,y)={}_{l}  \partial _{y} F(x,y), \\ 
& F(x,\, 0)=x^{\;n} , \\ 
& F(0,y)=y^{\;n}  \\ 
& {}_{l} \partial _{\xi } =\partial _{\xi }\; \xi \, \partial _{\xi }  \end{split}\end{equation} 
where ${}_{l} \partial _{\xi } $ is the Laguerre $l\!-\!derivative$, introduced in the past to deal with the monomiality properties of Laguerre polynomials \cite{D.Babusci}.\\

\noindent It is easily checked that the solution of eq. \eqref{GrindEQ__1_}, when $x,y>0$, can be cast in the form
 
\begin{equation} \label{GrindEQ__2_} 
\lambda _{n} (x,y)=\sum _{r=0}^{n}\binom{n}{r}^{2}  x^{\;n-r} y^{\;r}  
\end{equation} 
where $\lambda _{n} (x,y)$ is an example of hybrid polynomial, introduced in ref. \cite{G.Dattoli}.\\

\noindent For reasons which will be clear in the following, we introduce the notation, borrowed from ref. \cite{D.E.Edmunds},
\begin{equation} \label{GrindEQ__3_} 
\lambda _{n} (x,y)=\sum _{r=0}^{n}\binom{n}{r}^{2}  x^{\;n-r} y^{\;r} =(x\oplus _{l} y)^{\;n}  
\end{equation} 
to define a composition rule, which will be defined Laguerre binomial sum.\\

\noindent  It is evident that such a notion is an extension of the Newton Binomial, which can be generated by the action of the shifting exponential operator on an ordinary monomial, namely

\begin{equation}
e^{\;y\, \partial _{x} } \, x^{\;n} =\sum _{r=0}^{\infty }\dfrac{y^{\;r} }{r!} \partial _{x}^{\;r}  x^{\;n} =\sum _{r=0}^{n}\dfrac{y^{\;r} }{r!} \dfrac{n!}{(n-r)!}  x^{\;n-r} =(x+y)^{\;n}
\end{equation}            
An analogous rule for the generation of lbs can be achieved by replacing the exponential function with the $ l\!-\!exponential $

\begin{equation} \label{GrindEQ__5_} 
_{l} e(\eta )=\sum _{r=0}^{\infty }\frac{\eta ^{\;r} }{\left(r!\right)^{2} }   
\end{equation} 
and the ordinary derivative with the $l$-derivative, satisfying the identity

\begin{equation} \label{GrindEQ__6_} 
{}_{l} \partial _{\eta }^{\;n} =\partial _{\eta }^{\;n} \eta ^{\;n} \, \partial _{\eta }^{\;n}  
\end{equation} 
accordingly we find

\begin{equation}\begin{split} \label{GrindEQ__7_} 
& {}_{l} e(y\, {}_{l} \partial _{x} )\, x^{n} =\sum _{r=0}^{\infty }\frac{y^{r} }{\left(r!\right)^{2} } \partial _{x}^{r}  x^{r} \partial _{x}^{r} \, x^{n} = \\ 
& =\sum _{r=0}^{n}\frac{y^{r} }{\left(r!\right)^{2} } \frac{\left(n!\right)^{2} }{\left[(n-r)!\right]^{2} }  x^{n-r} =(x\oplus _{l} y)^{n} 
\end{split}\end{equation} 
The function ${}_{l} e(\eta )$ is a $0$-order Bessel Tricomi function \cite{Tricomi}\footnote{It can be expressed in terms of the $0$-order modified Bessel through the identity   ${}_{l} e(\eta )=I_{0} (2\, \sqrt{\eta } )$.} and satisfies the "Laguerre"-eigenvalue equation

\begin{equation} \label{GrindEQ__9_} 
\, {}_{l} \partial _{x} \left[{}_{l} e(\lambda \, x)\right]=\lambda \left[{}_{l} e(\lambda \, x)\right] 
\end{equation} 
which corroborates the interpretation of its role as that of an $l$-exponential.\\

\noindent According to the previous identities we can also state that

\begin{equation}\begin{split} \label{GrindEQ__9b_} 
& {}_{l}e(y\, {}_{l} \partial _{x} )\, {}_{l}e(x)={}_{l}e(y)\, {}_{l}e(x), \\ 
& {}_{l}e(y\, {}_{l} \partial _{x} )\, {}_{l}e(x)={}_{l}e(x\oplus _{l} y)
\end{split}\end{equation} 
allowing the derivation of the following "semi-group" property of the $l$-exponential

\begin{equation}
{}_{l} e(y)\, {}_{l} e(x)={}_{l} e(x\oplus _{l} y)                      
\end{equation} 

In full analogy with the ordinary Euler formula we introduce the l-trigonometric $l\!-\!t$ functions through the identity

\begin{equation} \label{GrindEQ__11_} 
\, {}_{l} e(i\, x)={}_{l} c(x)+i\, {}_{l} s(x) 
\end{equation} 
where $l\!-\!t$ \textit{cosine} and \textit{sine} functions are specified by the series\footnote{ The $l\!-\!h$ functions are defined by the corresponding series expansion   $\begin{array}{l} {{}_{l} ch(x)=\sum _{r=0}^{\infty }\frac{x^{2r} }{\left[(2\, r)!\right]^{2} }  ,} \\ {{}_{l} sh(x)=\sum _{r=0}^{\infty }\frac{x^{2r+1} }{\left[(2\, r+1)!\right]^{2} }  } \end{array}$. }

\begin{equation}\begin{split}\label{GrindEQ__12_}
& {}_{l} c(x)=\sum _{r=0}^{\infty }\frac{(-1)^{r} x^{2r} }{(2\, r)!^{2} }  , \\ 
& {}_{l} s(x)=\sum _{r=0}^{\infty }\frac{(-1)^{r} x^{2r+1} }{(2\, r+1)!^{2} }  
\end{split}\end{equation}                    
It is easily checked that they satisfy the identities

\begin{equation}\begin{split} \label{GrindEQ__13_} 
& {}_{l} \partial _{x} \left[{}_{l} c(\alpha \, x)\right]=-\alpha \, {}_{l} s(\alpha \, x), \\ 
& {}_{l} \partial _{x} \left[{}_{l} s(\alpha \, x)\right]=\alpha \, {}_{l} c(\alpha \, x)
 \end{split} \end{equation} 
and therefore the  "harmonic" equation

\begin{equation}\begin{split}
& \left({}_{l} \partial _{x} \right)^{2} \left[{}_{l} c(\alpha \, x)\right]=-\alpha ^{2} {}_{l} c(\alpha \, x), \\ 
& \left({}_{l} \partial _{x} \right)^{2} \left[{}_{l} s(\alpha \, x)\right]=-\alpha ^{2} {}_{l} s(\alpha \, x)
\end{split}\end{equation}

It is worth noting that Laguerre  and ordinary trigonometric functions are linked by the Borel like transforms \cite{Borel}

\begin{equation}\begin{split} \label{GrindEQ__70b_} 
& \int _{0}^{\infty }e^{-t}  {}_{l} c(xt)dt=\cos (x) \\ 
& \int _{0}^{\infty }e^{-t}  {}_{l} s(xt)dt=\sin (x)
\end{split} \end{equation} 
The use of the dilatation operator identity \cite{D.Babusci}

\begin{equation}
t^{x\, \partial _{x} } f(x)=f\left(t\, x\right)
\end{equation}
allows to cast transforms of the type \eqref{GrindEQ__70b_} in the operatorial form

\begin{equation}\begin{split}
& \int _{0}^{\infty }e^{-t}  f(xt)dt=\int _{0}^{\infty }e^{-t}  t^{x\, \partial _{x} } dt\, f(x)= \\ 
& =\Gamma (x\, \partial _{x} +1)\, f(x) 
\end{split}\end{equation}         
thus yielding the following identifications

\begin{equation}\begin{split} \label{GrindEQ__14_} 
& {}_{l} c(x)=\left(\Gamma (x\, \partial _{x} +1)\right)^{-1} \cos (x) \\ 
& {}_{l} s(x)=\left(\Gamma (x\, \partial _{x} +1)\right)^{-1} \sin (x)
\end{split} \end{equation} 

We note that

\begin{equation}\begin{split}
& {}_{l} e(i\, x)={}_{l}c(x)+i\, {}_{l}s(x)= \\ 
& =\left(\Gamma (x\, \partial _{x} +1)\right)^{-1} e^{i\, x}
\end{split} \end{equation} 	
 
The use of the properties of the $l$-exponential function allows the  derivation of the following addition theorems for the functions in eqs. \eqref{GrindEQ__12_} 

\begin{equation}\begin{split} \label{GrindEQ__15_} 
& {}_{l} c(x\oplus_{l} y)={}_{l} c(x){}_{l} c(y)-{}_{l} s(x)\, {}_{l} s(y), \\ 
& {}_{l} s(x\oplus_{l} y)={}_{l} c(x){}_{l} s(y)+{}_{l} s(x)\, {}_{l} c(y)
\end{split}\end{equation} 
The proof of the last identity is given below, by noting that

\begin{equation}\begin{split} \label{GrindEQ__16_} 
& {}_{l} e(i\, x)\, e(i\, y)=\left[{}_{l} c(x)+i\, {}_{l} s(x)\right]\, \left[{}_{l} c(y)+i\, {}_{l} s(y)\right]= \\ 
& =\left[{}_{l} c(x){}_{l} c(y)-{}_{l} s(x)\, {}_{l} s(y)\right]+i\, \left[{}_{l} c(x){}_{l} s(y)+{}_{l} s(x)\, {}_{l} c(y)\right]
\end{split}\end{equation} 
and since

\begin{equation} \label{GrindEQ__17_} 
{}_{l} e(i\, x)\, {}_{l} e(i\, y)={}_{l} e(i\left(x\oplus_{l} y\right))={}_{l} c(x\oplus_{l} y)+i\, {}_{l} s(x\oplus_{l} y) 
\end{equation} 
we can equate real and imaginary parts to infer the identities \eqref{GrindEQ__15_}.\\

 It is evident that, according to the procedure we have proposed, the properties of ordinary trigonometric functions are extended to their $l$-counterparts, provided that the ordinary sum is replaced by the composition rule specified in eq. \eqref{GrindEQ__3_}.\\

\noindent  The formalism allows the derivation of the corresponding duplication formulae, which can be stated by defining the following product rule 

 \begin{equation}
 \left( x\oplus_{l} x\right) ^{n}=\dfrac{(2n)!}{(n!)^{2}}x^{n}:=\left(2\otimes_{l} x \right) ^{n}
 \end{equation}                      
which along with eq. \eqref{GrindEQ__15_} yields 

\begin{equation}\begin{split}
& {}_{l} c(2\otimes_{l} x)=\left({}_{l} c(x)\right)^{2} -\left({}_{l} s(x)\right)^{2} , \\
 & {}_{l} s(2\otimes_{l} x)=2\;{}_{l} c(x){}_{l} s(x)
\end{split}\end{equation}
 Furthermore the sum can be iterated as

\begin{equation}\begin{split}
 & \left( x\oplus_{l} (x\oplus_{l} x)\right) ^{n}= \left( 3\otimes_{l} x\right) ^{n}, \\
& \left( x\oplus_{l} (..._{l}\; (x\oplus_{l} x))\right) ^{k}=\left( n\otimes_{l} x\right) ^{k} \end{split}\end{equation}               and, accordingly, we can state the following extension of the De Moivre formula

\begin{equation}
\left[{}_{l} c(x)+i\, {}_{l} s(x)\right]^{n} ={}_{l} c\left(n\otimes_{l} x\right)+i\, {}_{l} s\left(n\otimes_{l} x\right)
\end{equation}

 The Lissajous plot of the $ l\!-\!t $ functions is reported below and will be further commented in the concluding remarks.\\

%\begin{figure}[htp]
%	\centering
%	\begin{minipage}[b]{0.8\columnwidth}
%		\includegraphics[width=.9\textwidth]{FigLiss12}
%		\caption{\large   $ {}_{l} s(x) $ vs    $ {}_{l} c(x) $.}
%		\label{FigAss1}
%	\end{minipage}
%	\newline
%	\begin{minipage}[b]{0.78\columnwidth}
%		\includegraphics[width=.9\textwidth]{FigLiss3_4}
%		\caption{\large $ {}_{l} s(x\oplus_{l} y) $ vs $ {}_{l} c(x\oplus_{l} y) $.}
%		\label{FigAssC2}
%	\end{minipage}
%	\caption{\large Fish-like Lissajaux diagram of l-t functions for different parameters.}\label{fig1}
%\end{figure} 
%\begin{figure}[htp]
%\includegraphics[width=.9\textwidth]{FigLiss12}
%\caption{ Fish-like Lissajous diagram of $ l\!-\!t $ functions, $ {}_{l} s(x) $ vs    $ {}_{l} c(x) $, for different parameters.}\label{fig1}
%\end{figure} 
 
 The composition rule \eqref{GrindEQ__3_} is one of the pivotal points of this paper and will be further exploited in the following. It is therefore worth stressing some of its properties reported below

\begin{equation}\label{trivial}
\begin{split}
& a) \;\; (x\oplus _{l} y)^{n} =(y\oplus _{l} x)^{n} \\                  
& b) \;\; (1\oplus _{l} 1)^{n} =\dfrac{(2\, n)!}{\left(n!\right)^{2} } \\                  
& c) \;\; (1\oplus _{l} (-1))^{n}=\dfrac{i^{\;n}\;n! }{\left( \left(\frac{n}{2}   \right)!\right) ^{2} }\dfrac{\left( 1+(-1)^{n}\right)}{2}= \left\lbrace \begin{array}{ll} 0,&  n =2k+1,\; k\in \varmathbb{N}\\[1.1ex]  \dfrac{i^{\;n}\;n! }{\left( \left(\frac{n}{2}   \right)!\right) ^{2} }, & n=2k,\;\;\; k\in \varmathbb{N}  \end{array}\right. \\
& d) \;\; (i\oplus _{l} (-i))^{n}=\dfrac{(-1)^{n}\;n! }{\left( \left(\frac{n}{2}   \right)!\right) ^{2} }\dfrac{\left( 1+(-1)^{n}\right)}{2}= \left\lbrace \begin{array}{ll} =0,&  n =2k+1,\; k\in \varmathbb{N}\\[1.1ex] = \dfrac{(-1)^{n}\;n! }{\left( \left(\frac{n}{2}   \right)!\right) ^{2} }, & n=2k,\;\;\; k\in \varmathbb{N}  \end{array}\right.
\end{split}
\end{equation}
It is furthermore worth noting that 

\begin{equation}\label{Nep}
e) \;\;\;\; {}_{l}e( x)=\lim_{n\rightarrow\infty}\left(1\oplus _{l} \left( \dfrac{x}{n^{2}}\right)\right)^{n}   
\end{equation}
and that the quantity $\left(1\oplus _{l} \left( \dfrac{1}{n^{2}}\right)\right)^{n}$, a kind of Napier-Laguerre number, is presumably trascendent.\\

 The previous identities, albeit trivial, are important to appreciate the structural differences with respect to their ordinary counterparts, for example eq. \eqref{trivial} implies that
 
\begin{equation} \label{GrindEQ__21_} 
{}_{l}e(i\, x)\, {}_{l} e(-i\, x)\ne 1 
\end{equation} 
and therefore that

\begin{equation}
\left[{}_{l} c(x)\right]^{2} +\left[{}_{l} s(x)\right]^{2} \ne 1  
\end{equation} 
                    
In the forthcoming sections we will go deeper into the theory of these families of functions and we will be able to better appreciate the similitudes and the differences with the ordinary forms.

\begin{figure}[htp]
	\centering
	\includegraphics[width=.5\textwidth]{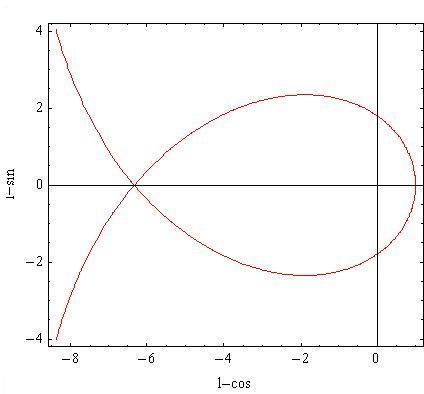}
	\caption{ Fish-like Lissajous diagram of $ l\!-\!t $ functions, $ {}_{l} s(x) $ vs    $ {}_{l} c(x) $.}\label{fig1}
\end{figure}

\section{ The generalized Trigonometric functions, ordinary and higher order Bessel functions}
\large

 The Bessel functions are characterized by a continuous variable and by a real or complex index. A fairly natural extension of the function defined by eqs. \eqref{GrindEQ__12_} is therefore provided by the $l\!-\!t$ functions, associated with the $\alpha$-order Tricomi- Bessel function \cite{Tricomi}
 
\begin{equation} \label{GrindEQ__23_} 
_{l} e_{\alpha } (\eta )=\sum _{r=0}^{\infty }\frac{\eta ^{\; r} }{r!\, \Gamma (r+\alpha +1)}   
\end{equation} 
which is an eigenfunction of the operator

\begin{equation}\begin{split} & {}_{(\alpha ,\, l)} \partial _{x} =\partial _{x} x\, \partial _{x} +\alpha \, \partial _{x} = \\
& =x^{-\alpha } \partial _{x} x^{\;\alpha +1} \partial _{x} 
\end{split}\end{equation}    
                 
We can now proceed as in the previous section, by noting that the polynomials
\begin{equation}\begin{split} \label{GrindEQ__25_} 
&\lambda _{n} (x,y;\alpha )=\sum _{r=0}^{n}\binom{n}{r}\; \dfrac{\Gamma (n+\alpha +1)}{\Gamma (n-r+\alpha +1)\, \Gamma (r+\alpha +1)\, }  x^{\;n-r} y^{\;r} = \\ 
&=(x\oplus _{(\alpha ,\, l)} y)^{\;n} 
 \end{split} \end{equation} 
are solutions of the equation

\begin{equation}
{}_{(\alpha ,\, l)} \partial _{x} l_{n} (x,y;\alpha )={}_{(\alpha ,\, l)} \partial _{y} l_{n} (x,y;\alpha ) 
\end{equation}                  
We further define the composition rule

\begin{equation}\label{comprule}
	{}_{l} e_{\alpha } (y\, {}_{\left(\alpha ,\, l\right)} \partial _{x} )\, x^{n} =(x\oplus _{(\alpha ,\;l)} y)^{n}
\end{equation}                
and prove that

\begin{equation}\begin{split}\label{GrindEQ__30b_} 
 &{}_{l} e_{\alpha } (y\, {}_{\left(\alpha ,\, l\right)} \partial _{x} )\, {}_{l} e_{\alpha } (x)=\, {}_{l} e_{\alpha } (x\oplus _{(\alpha ,\;l)} y), \\ 
& {}_{l} e_{\alpha } (y\, )\, {}_{l} e_{\alpha } (x)=\, {}_{l} e_{\alpha } (x\oplus _{(\alpha ,\;l)} y)
\end{split}\end{equation}             
Finally, by a slight extension of the discussion of the introductory section, we define the $l\!-\!t$ functions of $\alpha$-order as

\begin{equation}\begin{split}\label{GrindEQ__29_} 
&{}_{\left(\alpha ,\, l\right)} c(x)=\sum _{r=0}^{\infty }\dfrac{(-1)^{r} x^{2r} }{\left[(2\, r)!\right]\Gamma (\alpha +2r+1)}  , \\ 
& {}_{(\alpha ,\, l)} s(x)=\sum _{r=0}^{\infty }\dfrac{(-1)^{r} x^{2r+1} }{\left[(2\, r+1)!\right]\Gamma (\alpha +2r+2)} 
\end{split}\end{equation}
which are shown to satisfy the addition theorems \eqref{GrindEQ__15_} for the composition rule \eqref{comprule} .\\

 The procedure can be further extended by the use of Humbert Bessel like functions, which are defined by the series \cite{Cocolicchio}
\begin{equation} \label{GrindEQ__32b_} 
_{l} e_{\alpha ,\;\beta } (\eta )=\sum _{r=0}^{\infty }\frac{\eta ^{r} }{r!\, \Gamma (r+\alpha +1)\, \Gamma (r+\beta +1)}   
\end{equation} 
They satisfy the differential equation

\begin{equation} \label{GrindEQ__31_} 
\left[ \partial _{\eta } \left(\alpha +\eta \, \partial _{\eta } \right)\left(\beta +\eta \, \partial _{\eta } \right)\right]  {}_{l} e_{\alpha ,\beta } (\eta )= {}_{l} e_{\alpha ,\;\beta } (\eta ) 
\end{equation} 
and are therefore eigenfunctions of the operator

\begin{equation} \label{GrindEQ__32_} 
{}_{\left(\alpha ,\;\beta ,\, l\right)} \partial _{\eta } =\partial _{\eta } \eta \, \partial _{\eta } \eta \, \partial _{\eta } +\left(\alpha +\beta \right)\, {}_{l} \partial _{\eta } +\alpha \, \beta \, \partial _{\eta }  
\end{equation} 
whose  $ r$-order  derivative is

\begin{equation}\label{rder}
{}_{\left(\alpha ,\;\beta ,\, l\right)} \partial _{\eta }^{r}x^{n}=\dfrac{n!(n+\alpha)!(n+\beta)!}{(n-r)!(n-r+\alpha)!(n-r+\beta)!}x^{n-r}
\end{equation}
By following the same procedure as before, we introduce the composition rule

\begin{equation}\begin{split}\label{GrindEQ__33_} 
& \left( x\oplus _{(\alpha ,\, \beta ,\;l)} y\right) ^{n} = \\ 
& =\sum _{r=0}^{n}\binom{n}{r}\, \dfrac{\Gamma (n+\alpha +1)\, \Gamma (n+\beta +1)}{\Gamma (n-r+\alpha +1)\, \Gamma (r+\alpha +1)\, \Gamma (n-r+\beta +1)\, \Gamma (r+\beta +1)\, }  x^{\;n-r} y^{\;r} 
\end{split} \end{equation}  
so that the associated $l\!-\!t$ functions, defined as,

\begin{equation}\begin{split} \label{GrindEQ__34_} 
& {}_{\left(\alpha ,\;\beta ,\, l\right)} c(x)=\sum _{r=0}^{\infty }\dfrac{(-1)^{\;r} x^{\;2r} }{(2\, r)!\Gamma (\alpha +2r+1)\, \Gamma (\beta +2\, r+1)}  , \\ 
& {}_{(\alpha ,\;\beta ,\, l)} s(x)=\sum _{r=0}^{\infty }\dfrac{(-1)^{\;r} x^{\;2r+1} }{(2\, r+1)!\Gamma (\alpha +2r+2)\, \Gamma (\beta +2\, r+2)}
\end{split}\end{equation} 
are straightforwardly shown to satisfy the differential equations

\begin{equation}\begin{split} \label{GrindEQ__35_} 
& {}_{\left(\alpha ,\;\beta ,\, l\right)} \partial _{\eta } \left[{}_{\left(\alpha ,\;\beta ,\, l\right)} c(\lambda \, x)\right]=-\lambda \left[{}_{\left(\alpha ,\;\beta ,\, l\right)} s(\lambda \, x)\right]\, , \\ 
& {}_{\left(\alpha ,\;\beta ,\, l\right)} \partial _{\eta } \left[{}_{\left(\alpha ,\;\beta ,\, l\right)} s(\lambda \, x)\right]=\lambda \left[{}_{\left(\alpha ,\;\beta ,\, l\right)} c(\lambda \, x)\right]
\end{split}\end{equation} 
and the addition theorems based on the extension of the definition of sum specified in eq. \eqref{GrindEQ__9b_} and \eqref{GrindEQ__30b_}.\\

  We have shown that the concept of $ l\!-\!t $ function is a fairly natural consequence of the  notion of Laguerre derivative, of its extensions and of the associated eigenfunctions, which belong to Bessel like forms. In the following we will show how to frame the Cholewinsky-Reneke $ l\!-\!h $ functions within the present framework. Before entering this specific aspect of the problem we introduce some consequences of the previous formalism on the theory of diffusion equation associated to the Laguerre derivative and to its generalization.

\section{ Bessel Diffusion equations}

\noindent This short section, in which we will discuss some evolutive equations based on the operators introduced in the previous sections, is an apparent detour from the main stream of the paper .\\

 Laguerre type diffusive equations like \cite{D.Babusci}

\begin{equation}\begin{split}\label{GrindEQ__36_} 
&\partial _{\tau } F(x,\tau )={}_{l} \partial _{x} F(x,y), \\ 
& F(x,\, 0)=f(x)
 \end{split} \end{equation} 
can be formally solved as
\begin{equation} \label{GrindEQ__37_} 
F(x,\, \tau )=e^{\tau \, {}_{l} \partial _{x} } f(x) 
\end{equation} 
To make the above solution meaningful it is necessary to specify how to calculate the action of the exponential operator containing the Laguerre derivative on the function $f(x)$. We will discuss therefore, as introductory example, the case in which $f(x)=e^{x} $, and proceed as it follows:

\begin{enumerate}
\item  We note that the exponential can be written as an integral transform of the Tricomi function

\begin{equation} \label{GrindEQ__40b_} 
e^{x} =\int _{0}^{\infty }e^{-t}  {}_{l} e(x\, t)\, dt 
\end{equation} 

\item  We use the properties \eqref{GrindEQ__9b_} to end up with
\end{enumerate}
 
\begin{equation} \label{GrindEQ__41b_} 
e^{\tau \, {}_{l} \partial _{x} } e^{x} =\int _{0}^{\infty }e^{-t\;(1-\tau )}  {}_{l}e(x\, t)\, dt=\dfrac{1}{1-\tau } e^{\frac{x\, }{1-\tau } }  
\end{equation} 
More in general, whenever
 
\begin{equation} \label{GrindEQ__40_} 
f(x)=\int _{0}^{\infty }\tilde{f}(t) {}_{l} e(x\, t)\, dt 
\end{equation} 
the solution of the problem \eqref{GrindEQ__36_} can be cast in the form

\begin{equation} \label{GrindEQ__43b_} 
F(x,\, \tau )=\int _{0}^{\infty }\tilde{f}(t)\, e^{t\;\tau }  {}_{l} e(x\, t)\, dt 
\end{equation} 
and  $\tilde{f}(t)$ is the $l$-transform of the function $f(x)$.\\

 Before going further with the above formalism, we note that the equation

\begin{equation}\begin{split}\label{GrindEQ__44b_} 
&\partial _{\tau } F(x,\tau )={}_{(\alpha,\;\beta,\;l)} \partial _{x} F(x,y), \\ 
& F(x,\, 0)=f(x)
\end{split} \end{equation} 
can be solved in an analogous way provided that we replace ${}_{l}e(x)$ with ${}_{l}e_{\alpha ,\, \beta } (x\, )$ in eq.
% \eqref{GrindEQ__40b_} and
 \eqref{GrindEQ__43b_}.\\

\noindent In the case in which we consider equations of the type

\begin{equation}\begin{split}\label{GrindEQ__43_} 
& {}_{l} \partial _{\tau } F(x,\tau )={}_{l} \partial _{x} F(x,y), \\ 
& F(x,\, 0)=f(x)
\end{split} \end{equation} 	
the solution of the problem can be obtained as

\begin{equation}
F(x,\, \tau )= {}_{l}e(\tau \, {}_{l} \partial _{x} )f(x)=f(x\oplus_{l} \tau )
\end{equation}  
      Further comments on the previous statements will be provided in the forthcoming parts of the paper.

\section{ Pseudo-hyperbolic functions and generalized Airy Diffusion equations}

\noindent As already stressed the study of generalized forms of trigonometric and of hyperbolic functions is an old ``leit motiv'' in the mathematical literature. On the eve of the seventies of the last century Ricci introduced \cite{P.E.Ricci} a family of pseudo hyperbolic functions (PHF), which will be proven of noticeable importance for the topics we are discussing. \\

\noindent According to ref. \cite{P.E.Ricci} the PHF of order $ 3 $ are defined by the series
 
\begin{equation}\begin{split}
& {}_{[k,3]} e\left(x\right)=\sum _{r=0}^{\infty }\dfrac{x^{3\, r+k} }{(3\, r+k)!}  , \\ 
& k=0,1, 2 
\end{split}\end{equation}                          
On account of the properties of the cubic roots of the unit
 
\begin{equation}\begin{split} \label{GrindEQ__46_} 
& \hat{\omega }_{p} =e^{\frac{2\, i\, p\;\pi }{3} \, } ,\;\;\;\; p=0,1,2, \\ 
& \hat{\omega }_{p}^{3} =1,\; \;\;\;\;\;\;\;\;p=0, 1, 2,\\
& \hat{\omega }_{p}^{2} +\hat{\omega }_{p}=-1, \;\;\;p=1,2
\end{split}\end{equation} 
 the following Euler-like exponential formulae can be stated \cite{Migliorati}

\begin{equation}\begin{split} 
& e^{\hat{\omega}\, x} =\sum _{k=0}^{2}\hat{\omega}^{k}\; {}_{[k,3]} e (x), \\ 
& {}_{[k,3]}e(x)=\dfrac{1}{3} \sum _{p=0}^{2}\hat{\omega }_{p}^{k}\; e^{\hat{\omega }_{p} x}
\end{split}\end{equation}  
 
The PHF of order 3 are eingenfunctions of the cubic operator

\begin{equation} \label{GrindEQ__48_} 
\left(\partial {}_{x} \right)^{3} {}_{[k,3]} e(\lambda \, x)=\lambda ^{3} {}_{[k,3]} e(\lambda \, x) 
\end{equation} 
and can be exploited to generalize the exponential translation operator as

\begin{equation}\label{51b}
{}_{[k,3]} \hat{T}(y)={}_{[k,3]} e(y\, \partial _{x} ) 
\end{equation}                  
In the case of  $k=0$ the action of this operator on an ordinary monomial is given by 

\begin{equation}\begin{split}\label{52b}
& {}_{[0,3]} \hat{T}(y)\, x^{3n} =\frac{1}{3} \left[e^{\hat{\omega }_{0} \, y\, \partial _{x} } +e^{\hat{\omega }_{1} \, y\, \partial _{x} } +e^{\hat{\omega }_{2} \, y\, \partial _{x} } \right]\, x^{3n} = \\ 
& =\frac{1}{3} \sum _{\alpha =0}^{2}(x+\hat{\omega }_{\alpha }  y)^{3n} =(x\oplus_{\left[0,3\right]} y)^{3n} 
 \end{split}\end{equation}   
By direct use of the series expansion definition of the function $_{[0,3]} e\left(x\right)$ , we end up with\footnote{A straightforward consequence of eq. \eqref{52b} is that (see also ref. \cite{Cholewinski})\\
$(1\oplus_{[0,3]} 1)^{3n} =\dfrac{1}{3} \left(2^{3n}+\left(1+e^{\frac{2i\pi}{3}} \right)^{3n}+\left(1+e^{\frac{4i\pi}{3}} \right)^{3n}  \right) =\dfrac{1}{3} \left(2^{3\, n} +(-1)^{n} 2\right)$ }

\begin{equation}\begin{split}\label{53b}
& {}_{[0,3]} e\left(y\, \partial _{x} \right)x^{3\, n} =\sum _{r=0}^{\infty }\frac{y^{3r} }{(3r)!} \,  \partial _{x}^{3r} x^{3n} = \\ 
& =\sum _{r=0}^{n}\binom{3n}{3r} \; y^{3r}  x^{3\, \left(n-r\right)} =(x\oplus_{\left[0,3\right]} y)^{3n} 
\end{split}\end{equation}                                        
It is therefore evident that the following further  identities can be stated
\begin{equation}\begin{split} \label{GrindEQ__52_} 
& {}_{[0,3]} e\left(y\, \partial _{x} \right)\, {}_{[0,3]} e(x)={}_{[0,3]} e(x\oplus_{0} y), \\ 
& {}_{[0,3]} e\left(y\, \partial _{x} \right)\, {}_{[0,3]} e(x)={}_{[0,3]} e(y)\, {}_{[0,3]} e(x)  
\end{split}\end{equation} 
which, once merged, yield 
\begin{equation} \label{GrindEQ__53_} 
{}_{[0,3]} e\left(x\right)\, {}_{[0,3]} e(y)={}_{[0,3]} e(x\oplus_{[0,3]} y) 
\end{equation} 
In this way we have obtained a result allowing the introduction of $t\!-\!h$ like functions according to the paradigm developed so far.\\

 By a straightforward generalization of the discussion developed in the introductory remarks we introduce the generalized  $h$-functions defined as
 
\begin{equation}\begin{split} \label{GrindEQ__54_} 
& {}_{[0,3]}ch\left(x\right)=\frac{1}{2} \left( {}_{[0,3]}e(x)+ {}_{[0,3]}e\left(-x\right)\right) =\sum _{r=0}^{\infty }\dfrac{x^{\;6\, r} }{(6\, r)!}   \\ 
& {}_{[0,3]}sh\left(x\right)=\frac{1}{2} \left( {}_{[0,3]} e(x)- {}_{[0,3]}e\left(-x\right)\right) =\sum _{r=0}^{\infty }\dfrac{x^{\;6\, r+3} }{(6\, r+3)!}  
\end{split}\end{equation} 
and easily state that the relevant addition theorems read

\begin{equation}\begin{split} \label{GrindEQ__55_} 
& {}_{[0,3]}ch\left( \alpha \oplus_{[0,3]} \beta \right) ={}_{[0,3]}ch(\alpha ){}_{[0,3]}ch(\beta )+ {}_{[0,3]} sh(\alpha ) {}_{[0,3]}sh(\beta ), \\
& {}_{[0,3]} sh\left( \alpha \oplus_{[0,3]} \beta \right) = {}_{[0,3]}ch(\alpha ) {}_{[0,3]}sh(\beta )+ {}_{[0,3]}sh(\alpha ) {}_{[0,3]}ch(\beta )  
\end{split}\end{equation}
furthermore, as discussed later in the paper, analogous  conclusions can be reached for ${}_{[k,m]} e(x),\,\; k=0,\dots, m-1$.\\

 To obtain the link with the topics discussed so far, we consider a particular case of the function defined in eq. \eqref{GrindEQ__32b_}, namely
 
\begin{equation} \label{GrindEQ__56_} 
_{l} e_{\alpha -\frac{1}{3} ,\, -\frac{2}{3} } \left(\left(\frac{\eta }{3} \right)^{3} \right)=\sum _{r=0}^{\infty }\frac{\eta ^{3r} }{3^{3r} r!\, \Gamma \left(r+\alpha +\frac{2}{3} \, \right)\Gamma \left(r+\frac{1}{3} \right)}   
\end{equation} 
The associated Laguerre derivative can be written in terms of the operator

\begin{equation} \label{GrindEQ__57_} 
{}_{\alpha } \hat{\vartheta }={}_{\left(\alpha -\frac{1}{3} ,-\frac{2}{3} ,\, l\right)} \partial _{\left(\frac{\eta }{3} \right)^{3} } =\partial _{\eta } \eta ^{-3\, \alpha } \partial _{\eta } \eta ^{3\, \alpha } \partial _{\eta }  
\end{equation} 
appearing in the generalized Airy equation

\begin{equation} \label{GrindEQ__58_} 
\partial _{t} F(x,t)={}_{\alpha } \hat{\vartheta }\, F(x,t) 
\end{equation} 
studied in ref. \cite{Cholewinski}.\\

 The function in eq. \eqref{GrindEQ__56_} is equivalent, apart from an unessential normalizing factor, to the function exploited by Cholewinsky and Reneke \cite{Cholewinski} to study the solution of eq. \eqref{GrindEQ__58_} which will be linked, in particular, to the following expression

\begin{equation}\label{61}
G_{\alpha } (\eta )=\Gamma 
\left(\frac{1}{3} \right)\, \Gamma \left(\alpha +\frac{2}{3} \right) {}_{l}e_{\alpha -\frac{1}{3} ,\, -\frac{2}{3} } \left(\left(\frac{\eta }{3} \right)^{3} \right) 
\end{equation}               
We can make the previous results more transparent by using the identity

\begin{equation} \label{GrindEQ__60_} 
\frac{1}{(3\, n)!} \Gamma \left(n+\frac{2}{3} \right)=\frac{\Gamma \left(\frac{1}{3} \right)\, \Gamma \left(\frac{2}{3} \right)}{3^{3n} n!\, \Gamma \left(n+\frac{1}{3} \right)}  
\end{equation} 
which allows to recast eq. \eqref{61} in the form

\begin{equation}\begin{split} \label{GrindEQ__61_} 
& G_{\alpha } (\eta )=\frac{1}{B\left(\frac{2}{3} ,\alpha \right)} \sum _{r=0}^{\infty }\frac{B\left(r+\frac{2}{3} ,\, \alpha \right)\, \eta ^{3r} }{\left(3r\right)\, !\, }  , \\ 
& B(x,y)=\int _{0}^{1}t^{x-1}  (1-t)^{y-1} dt
\end{split} \end{equation} 
and derive the identity

\begin{equation} \label{GrindEQ__66_} 
{}_{\alpha } \hat{\vartheta }\, G_{\alpha } (\lambda \;\eta )=\lambda ^{3} G_{\alpha } (\lambda \;\eta ) 
\end{equation}

Finally, the use of the Euler dilatation operator

\begin{equation} \label{GrindEQ__62_} 
t^{x\, \partial _{x} } f(x)=f(\, t\, x) 
\end{equation} 
 yields the following integral transform defining the function \eqref{61} in terms of the Pseudo hyperbolic function of order 3
 
\begin{equation}\label{65}
G_{\alpha } (\eta )=\frac{1}{B\left(\frac{2}{3} ,\alpha \right)} \int _{0}^{1}t^{-\frac{1}{3}} (1-t)^{\alpha -1 } {}_{[0,3]} e(\eta \, t^{\frac{1}{3} } )dt 
\end{equation}  
              
According to the paradigm developed so far we use the associated translation operator to define the composition rule 

\begin{equation}\begin{split} \label{GrindEQ__64_} 
& G_{\alpha } (y\; {}_{\alpha }\hat{\vartheta }^{\frac{1}{3}} )\, x^{3\, n} = \\ 
& =\sum _{r=0}^{n}\binom{n}{r}\, \frac{\Gamma \left(\alpha +\frac{2}{3} \right)\, \Gamma \left(\frac{1}{3} \right)\, \Gamma \left(n+\alpha +\frac{2}{3} \right)\Gamma \left(n+\frac{1}{3} \right) }{\Gamma \left(r+\frac{1}{3} \right)\, \Gamma \left(n-r+\frac{1}{3} \right)\, \Gamma \left(n-r+\alpha +\frac{2}{3} \right)\, \Gamma \left(r+\alpha +\frac{2}{3} \right)}  y^{3\, r} x^{3n-3 r} =\\
& =\left(x\oplus_{{}_{\alpha } [0|3]} y\right)^{3n}
\end{split} \end{equation} 
forming the generalized lbs for the addition theorem of the relevant $l\!-\!h$ functions.\\

%\noindent On account of the transform \eqref{65} we can also establish that
%\begin{equation} \label{GrindEQ__65_} 
%\left(x\oplus_{{}_{\alpha } [0|3]} y\right)^{3n} =\frac{1}{B\left(\frac{2}{3} ,\alpha \right)} \int _{0}^{1}t (1-t)^{\alpha -\frac{2}{3} } \left(x\oplus {}_{[0|3]} y\, t^{\frac{1}{3}} \right)^{3n} dt 
%\end{equation} 
A natural extension of the previously developed formalism allows the introduction of the l/h-functions; their definitions and properties, apart from some computational complications, does not produce any significant conceptual progress.\\

In this paper we have seen how a "wise" combination of different methods borrowed from special function theory, operational and umbral calculus and integral transforms, opens new interesting possibilities for the introduction and systematic study of new families of trigonometric functions.\\
In addition, the method yields, as byproduct, the opportunity of getting natural solutions of a large family of PDE belonging to the family of generalized heat equation, whose links with the radial heat equation have been touched on in ref. \cite{D.E.Edmunds} and will be more carefully discussed elsewhere.\\

Before closing the paper we like to stress two points just touched on in the previous part of the paper.\\
We have noted in eq. \eqref{Nep} that the Laguerre exponential can be obtained through a limit procedure analogous to that involving the Napier number. The same method can e.g. be used to state that

\begin{equation}
J_{0}(x)=\lim_{n\rightarrow\infty}\left(1\oplus_{l} \left(-\left( \dfrac{x}{2n} \right)^{2}\right)   \right)^{n} 
\end{equation}
which is recognized as the asymptotic limit of Laguerre polynomials \cite{L.C.Andrews}.\\

\noindent The second point we want to emphasize is the geometrical interpretation of the $ l\!-\!t $ functions from the geometrical point of view.\\
Such an interpretation is provided in Fig. \ref{fig1}, where we have considered the Lissajous curves plotting $l$-sin vs $l$-cos.\\

\begin{figure}[htp]
	\centering
	\includegraphics[width=.45\textwidth]{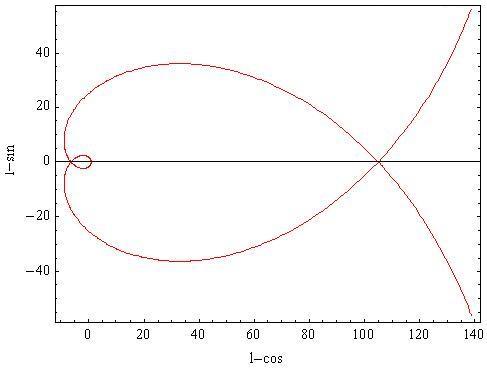}
	\caption{ Fish-like Lissajous diagram of $ l\!-\!t $ functions, $ {}_{l} s(x) $ vs    $ {}_{l} c(x) $ for larger $ x $-range.}\label{fig2}
\end{figure}

In Fig. \ref{fig2} it is evident that the curves are open since no periodic behaviour is envisaged. However a kind of self-similarity can be noted when the amplitude of the oscillations increase with increasing $x$.   \\

The last figure provides the explicit correspondence of the $l$-sinus and $l$-cosinus along with the relevant "$l$-angle", intended as the area intercepted by the $l$-curve and the segment forming the angle in the positive abscissa direction.

\begin{figure}[htp]
	\centering
	\includegraphics[width=.45\textwidth]{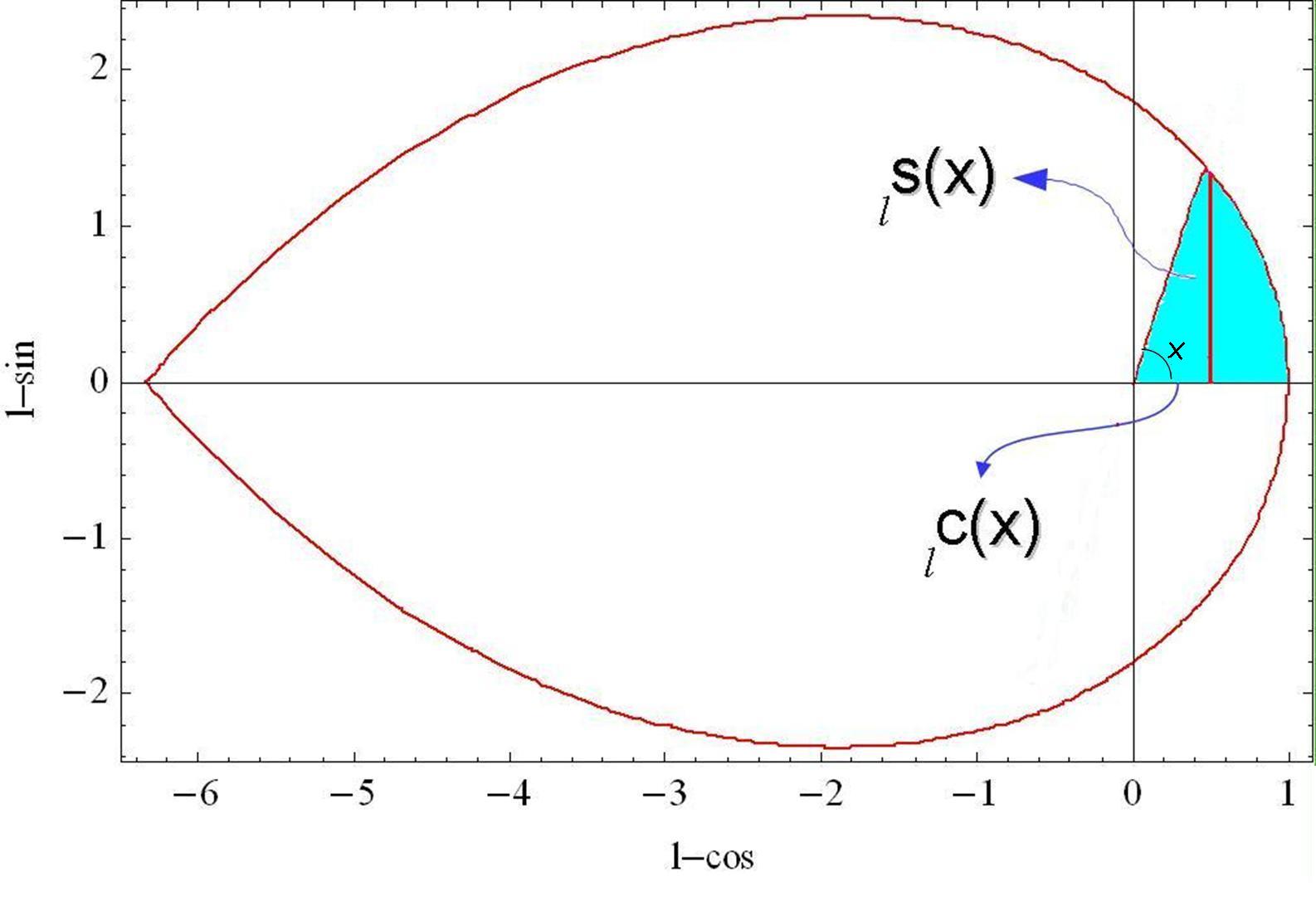}
	\caption{$l$-sinus and $l$-cosinus along with parameter $ x $ which is understood as the blue dashed area.}\label{fig3}
\end{figure} 
 \vspace{1.5cm}
 
 This paper has provided the formalism for a wider understanding of concepts associated with the Laguerre derivative, the underlying algebraic structures and  possible extensions to new forms of Trigonometry. In a forthcoming investigation we will provide an account of application of the method to problems in Physics and applied mathematics. \\


\textbf{References}

\begin{thebibliography}{}

\bibitem{L.C.Andrews} L.C. Andrews, \textit{Special functions for Engineers and Applied Mathematicians}, Mac Millan New York (1985).

\bibitem{E.Ferrari} E. Ferrari, Bollettino, UMI 18B, 933, (1981). 

\bibitem{D.E.Edmunds} D.E. Edmunds, P. Gurka, J. Lang,  \textit{Properties og generalized trigonometric functions}, Journal of Approximation Theory, 164, 47-56, (2012).

\bibitem{Gorska} D. Babusci, G. Dattoli, K. Gorska and K. Penson, \textit{"Symbolic Methods for the Derivation of Bessel Sum rules"}, J. Math. Phys. 54,  073501, 2013.

\bibitem{Cholewinski} F.M. Cholewinski, J.A.Reneke, \textit{"The Generalized Airy Diffusion Equation"}, Electronic Journal of Differential Equations, , No. 87, pp. 1-64, ISSN: 1072-6691, Vol. 2003(2003).

\bibitem{D.Babusci} D. Babusci,, G. Dattoli, M. Del Franco, \textit{Lectures on Mathematical Methods for Phisics}, RT/2010/58/ENEA.

\bibitem{G.Dattoli} G. Dattoli, S. Lorenzutta, P.E Ricci \& C. Cesarano, "On a family of hybrid polynomials", Integral Transforms and Special Finctions, 15:6, 485-490, DOI:10.1080/10652460412331270634, (2004).

\bibitem{Tricomi} F.G. Tricomi, "Funzioni Speciali", pp 408, Gheroni, (1959).

\bibitem{Borel} G. Dattoli, E. di Palma, E. Sabia, K. Górska, A. Horzela, K. A. Penson, \textit{"Operational Versus Umbral Methods and the Borel Transform"}, International Journal of Applied and Computational Mathematics, DOI 10.1007/s40819-017-0315-7, Feb 2017. 

\bibitem{Cocolicchio} D.P. Cocolicchio, G. Dattoli, H.M. Srivastava, "Proceedings of the Workshop Advanced Special Functions and Applications", Melfi School on Advanced Topics in Mathematics and Physics, Melfi (PZ), Italy, 9-12 May 1959, Aracne.

\bibitem{P.E.Ricci} P.E. Ricci, \textit{"Le funzioni pseudo-iperboliche e pseudo-trigonometriche"}, Istit. 
Mat. Appl., Fac.  Ingr. Univ.  Stud. Roma, Quad. 12, 37-49, (1978).

\bibitem{Migliorati} G.Dattoli, M.Migliorati, P.E.Ricci, \textit{"The Eisentein group and the pseudo hyperbolic functions"}, RT/2007/22/FIM.

\bibitem{K.Penson} D. Babusci, G. Dattoli, K. Gorska and K. Penson, "Lacunary Generating Functions for Laguerre Polynomials", arXiv:1302.4894v1 [math-ph], Journal of Combinatorial Analysi, to be published.

\end{thebibliography}
\end{document}